\documentclass[12pt,fleqn]{article}

\usepackage{geometry}

\usepackage[cp1251]{inputenc}
\usepackage[T2A]{fontenc}
\usepackage[english]{babel}

\usepackage{amsmath}
\usepackage{amssymb}
\usepackage{amsthm}
\usepackage{array}
\usepackage[auth-sc]{authblk}
\usepackage{bm}
\usepackage{cite}
\usepackage{color}
\usepackage{dsfont}
\usepackage{enumerate}
\usepackage{enumitem}
\usepackage{epsf}
\usepackage{euscript}
\usepackage{graphicx} 
\usepackage{indentfirst}
\usepackage{latexsym}
\usepackage{mathrsfs}
\usepackage{mathtools}
\usepackage{tikz}
\usepackage{pgf}

\newtheorem{theorem}             {Theorem}

\newlength{\templength}

\newcommand{\dis}{\vee}

\newcommand{\imp}{\rightarrow}

\begin{document}

\selectlanguage{english}

\title{Algorithmic properties of $\mathbf{QK4.3}$ and
  $\mathbf{QS4.3}$}

\author[1]{M.\,Rybakov}
\author[2]{D.\,Shkatov}
\affil[1]{IITP RAS, HSE University, Tver State University}
\affil[2]{University of the Witwatersrand, Johannesburg}

\date{} 

%
%

\maketitle

\begin{abstract}
  We prove that predicate modal logics $\mathbf{QK4.3}$ and
  $\mathbf{QS4.3}$ are undecidable---more precisely,
  $\Sigma^0_1$-complete--in languages with two individual variables,
  one modandic predicate letter, and one proposition letter.  
\end{abstract}


It is natural to expect that predicate modal logics should be
algorithmically harder than the classical predicate logic
$\mathbf{QCl}$, just as propositional modal logics are, as a rule,
algorithmically harder than the classical propositional logic.
Nevetheless, numerous predicate modal logics are just as hard as
$\mathbf{QCl}$, i.e., $\Sigma^0_1$\nobreakdash-complete: some---such
as $\mathbf{QK}$, $\mathbf{T}$, $\mathbf{S4}$, and $\mathbf{S5}$---are
recursively axiomatizable over $\mathbf{QCl}$~\cite{HC96,GShS}; others
are recursively embeddable~\cite{RShAiML18,RSh20JLC} into
$\mathbf{QCl}$ through the standard traslation.  It, however, turns
out that $\Sigma^0_1$\nobreakdash-complete predicate modal logics can
be distinguished from $\mathbf{QCl}$ by algorithmic properties of
their fragments: while the monadic fragment of $\mathbf{QCl}$ is
decidable~\cite{Lowenheim15,Behmann22}, the monadic fragments of most
$\Sigma^0_1$-complete modal logics are not~\cite{Kripke62}; while the
two-variable fragment of $\mathbf{QCl}$ is
decidable~\cite{Mortimer75,GKV97}, the two-variable fragments of most
$\Sigma^0_1$\nobreakdash-complete modal predicate logics are
not~\cite{KKZ05,RSh19SL}.  This leads to the study of the algorithmic
properties of the {\em fragments}\/ of modal predicate logics.

The algorithmic properties of one-variable and two-variable fragments
of first-order modal logics are also of interest due to close links
between those fragments and, respectively, two-dimensional and
three-dimensional propositional modal
logics~\cite{GSh98,GKWZ,Shehtman19,ShSh20,RShJLC21a}.

The study of the algorithmic properties of fragments of predicate
modal, and related superintuitionistic, logics
(see~\cite{Kripke62,MMO65,Mints68,Ono77,Gabbay81,AD90,GSh93,WZ01,KKZ05,RSh19SL,RShJLC20a,RShJLC21b};
for a summary of results, see~\cite[Introduction]{RShJLC20a}) is much
less advanced than similar research for $\mathbf{QCl}$~\cite{BGG97}.

In the present paper, we attempt to identify the minimal undecidable
fragments of the predicate counterparts $\mathbf{QK4.3}$ and
$\mathbf{QS4.3}$ of the well-know propositional modal logics
$\mathbf{K4.3}$ and $\mathbf{S4.3}$.  It is known~\cite{Corsi93} that
$\mathbf{QK4.3}$ and $\mathbf{QS4.3}$ are finitely axiomatizable over
$\mathbf{QCl}$ and, hence, they are $\Sigma^0_1$-complete.  The logics
$\mathbf{QK4.3}$ and $\mathbf{QS4.3}$ are faithfully characterized
using Kripke semantics with expanding domains~\cite{HC96}, \cite[\S
3.1]{GShS}: they are determined by all, respectively, strict and
partial linear orders~\cite{Corsi93}.  A closely related logic
$\mathbf{QK4.3.D.X}$ is determined by the rationals with the natural
strict order, viewed as a Kripke frame~\cite{Corsi93}.

The main interest of the results presented here is due the techniques
used: the known techniques for proving lower bounds in predicate modal
and superintuitionistic logics in languages with a few variables and a
few predicate letters~\cite{RSh19SL,RShJLC20a,RShJLC21b}, being based
on propositional-level
techniques~\cite{Halpern95,ChRyb03,Rybakov06,Rybakov08,RShICTAC18,RShIGPL18,RShIGPL19}
developed for logics of frames with unbounded branching, are
inapplicable to logics of linear frames (the only exception being our
earlier work~\cite{RSh20AiML}, where the techniques used for
establishing results reported here originate).

Predicate modal languages are obtained by enriching the classical
preciate language with a unary modal connective $\Box$ (for more
backgroud on predicate modal logic, see~\cite{HC96,FM98,GShS}

To recall the definitions of the logics we study, we use the following
notation: if $\Gamma$ is a set of formulas and $\varphi$ is a formula,
$\Gamma\oplus\varphi$ denotes the closure of $\Gamma\cup\{\varphi\}$
under modes ponens, generalization, necessitation, and predicate
substitution.  Then,

\settowidth{\templength}{\mbox{$\mathbf{QK4.3}$}}
$$
\begin{array}{lcl}
\mathbf{QK}
  & =
  & \parbox{\templength}{$\mathbf{QCl}$}
    ~\oplus~ \Box (p\to q) \imp (\Box p\to \Box q);
  \smallskip\\
\mathbf{QK4}
  & =
  & \parbox{\templength}{$\mathbf{QK}$}
    ~\oplus~ \Box p \imp \Box \Box p;
  \smallskip\\
\mathbf{QS4}
  & =
  & \parbox{\templength}{$\mathbf{QK4}$}
    ~\oplus~ \Box p \imp p;
  \smallskip\\
\mathbf{QS4.3}
  & =
  & \parbox{\templength}{$\mathbf{QS4}$}
    ~\oplus~ \Box (\Box p \imp q) \dis \Box (\Box q \imp p);
  \smallskip\\
\mathbf{QK4.3}
  & =
  & \parbox{\templength}{$\mathbf{QK4}$}
    ~\oplus~ \Box (p\wedge\Box p \imp q) \dis \Box (q\wedge\Box q \imp p);
  \smallskip\\
\mathbf{QK4.3.D.X}
  & =
  & \parbox{\templength}{$\mathbf{QK4.3}$}
    ~\oplus~ \Diamond \top
    ~\oplus~ \Box \Box p \imp \Box p.
\end{array}
$$



\begin{theorem}
  \label{thr:1}
  Logics $\mathbf{QK4.3}$, $\mathbf{QS4.3}$, and $\mathbf{QK4.3.D.X}$
  are $\Sigma^0_1$\nobreakdash-complete in languages containing one
  monadic predicate letter, one proposition letter, and two individual
  variables.
\end{theorem}

Theorem~\ref{thr:1} is proved as follows: we encode a
$\Sigma^0_1$\nobreakdash-hard tiling problem~\cite{Berger66} using
predicate modal formulas with only two variables, and only binary and
unary predicate letters; then, we simulate binary letters with monadic
one; finally, we simulate the monadic letters with just one monadic
and one nullary letter; the latter reductions do not use more than two
individual variables.

This result can be extended to logics containing the Barcan formula
$\mathit{bf} = \forall x\,\Box P(x)\to \Box\forall x\, P(x)$:

\begin{theorem}
  \label{thr:2}
  Logics $\mathbf{QK4.3}\oplus\mathit{bf}$,
  $\mathbf{QS4.3}\oplus\mathit{bf}$, and
  $\mathbf{QK4.3.D.X}\oplus\mathit{bf}$ are
  $\Sigma^0_1$\nobreakdash-complete in the language containing one
  monadic predicate letter, one proposition letter, and two individual
  variables.
\end{theorem}

Furthermore, we obtain the following generalisation of
Theorems~\ref{thr:1} and~\ref{thr:2}:

\begin{theorem}
  \label{thr:3}
  Every logic in the intervals
  $[\mathbf{QK4.3},\mathbf{QK4.3.D.X}\oplus\mathit{bf}]$ and
  $[\mathbf{QS4.3},\mathbf{QS4.3}\oplus\mathit{bf}]$ are
  $\Sigma^0_1$\nobreakdash-hard in the language containing one monadic
  predicate letter, one proposition letter, and two individual
  variables.
\end{theorem}

Our proofs of Theorems~\ref{thr:1}--\ref{thr:3} rely on Kripke
completeness of $\mathbf{QK4.3}$, $\mathbf{QS4.3}$, and
$\mathbf{QK4.3.D.X}$ and on soundness of
$\mathbf{QK4.3}\oplus\mathit{bf}$, $\mathbf{QS4.3}\oplus\mathit{bf}$,
and $\mathbf{QK4.3.D.X}\oplus\mathit{bf}$ with respect to,
respectively, $\mathbf{QK4.3}$-frames, $\mathbf{QS4.3}$-frames, and
$\mathbf{QK4.3.D.X}$-frames with constant domains.  

\bigskip

\textit{
  This research has been supported by the Russian Science Foundation
  with grant \mbox{21--18--00195}; it has been carried out at Tver
  State University.
}


\end{document}